\documentclass[11pt,reqno]{amsart}
\usepackage{amssymb,latexsym}
\usepackage{graphicx}
\usepackage{amsmath,amsthm}
\usepackage{enumerate}
\usepackage{fancyhdr,amssymb}
\usepackage{url}	
\usepackage{tikz}
\usepackage{pgf}

\usepackage{tipa}
\usepackage{stmaryrd}
\usepackage{amsmath,amsfonts}
\usepackage{mathtools} 
\usepackage{booktabs}

\newtheorem{thm}{Theorem}[section]
\newtheorem{cor}[thm]{Corollary}

\newtheorem{ques}[thm]{Question}

\newtheorem{rem}[thm]{Remark}

\newcommand{\Z}{{\mathbb Z}}

\newcommand{\Q}{{\mathbb Q}}

\newcommand{\Pc}{{\mathcal P}}

\newcommand{\Tc}{{\mathcal T}}

\newcommand{\ta}{{\theta}}
\newcommand{\ent}{E_{N}^\ta}

\begin{document}

\setcounter{page}{1}

\title[The $\ta$-congruent number elliptic curves via a Fermat-type theorem]{The $\ta$-congruent numbers elliptic curves via a Fermat-type theorem}
\author{Sajad Salami \and Arman Shamsi Zargar}
\address{Inst\'{i}tuto da Matem\'{a}tica e Estat\'{i}stica, Universidade Estadual do Rio de Janeiro (UERJ), Rio de Janeiro, Brazil}
\email{sajad.salami@ime.uerj.br}
\address{Department of Mathematics and Applications, Faculty of Basic Sciences, University of Mohaghegh Ardabili, Ardabil 56199-11367, Iran}
\email{zargar@uma.ac.ir}

\maketitle
\begin{abstract}
A positive integer $N$ is called a $\theta$-congruent number if there is a $\ta$-triangle
$(a,b,c)$ with rational sides for which the  angle between $a$ and $b$ is  equal to $\theta$ and its area is $N \sqrt{r^2-s^2}$, where $\theta \in (0, \pi)$,  $\cos(\theta)=s/r$, and $0 \leq |s|<r$ are coprime integers. It is attributed to Fujiwara \cite{fujw1} that
$N$ is a $\ta$-congruent number if and only if the elliptic curve $E_N^\ta: y^2=x (x+(r+s)N)(x-(r-s)N)$  has 	a point of order greater than $2$ in its group of rational points. Moreover, a natural number $N\neq 1,2,3,6$  is a $\ta$-congruent number if and only if rank of $E_N^\ta(\Q)$ is greater than zero.

In this paper, we answer positively to a question concerning with the existence of   methods  to create
new rational $\ta$-triangle for a $\ta$-congruent number $N$ from given ones by
generalizing the Fermat's algorithm, which produces new rational right triangles  for congruent numbers from a given one,
for any angle $\ta$ satisfying the above conditions.
We show that this generalization is analogous  to the duplication formula in   $E_N^\ta(\Q)$.
Then, based on the addition of two distinct points  in $E_N^\ta(\Q)$, we provide a way to find new rational $\ta$-triangles for the $\ta$-congruent number 	$N$ using  given two distinct ones.
Finally, we give an alternative proof for the Fujiwara's theorem~\ref{fuj1} and  one side of  Theorem~\ref{fuj2}.
In particular, we provide a list of all torsion points in $E_N^\ta(\Q)$ with corresponding  rational $\ta$-triangles.
\end{abstract}
\vspace{0.5cm}
{\it Keywords:} {$\ta$-congruent number, elliptic curve, Fermat type theorem,\\
Pythagorean $\ta$-triples}

{\it Subject class [2010]:} {Primary  11G05, Secondary 14H52}


\section{Introduction}
\label{intro}

A positive integer $N$ is called a \textit{congruent number}  if
it is equal to the area of a right triangle with rational sides;
equivalently, if there exist positive rational numbers $a$, $b$ and
$c$ such that $a< b <c,$ and
\begin{equation*}
a^2+b^2=c^2, \quad  ab=2N.
\end{equation*}
Determining all congruent numbers is an old problem in the theory of numbers and
there are different types of generalizations in the literature, see for example
\cite{fujw1,fujw2,kob1}. Specifically, in \cite{fujw1}, Fujiwara introduced the notion of $\theta$-congruent numbers  for an  angle $\ta \in (0, \pi)$ with rational cosine as a generalization of  the usual congruent numbers.
To be precise,  assume that $\ta \in (0, \pi)$
is an angle of a triangle given by  a triple $(a, b, c)$ of positive rational numbers for which  $\cos(\ta)=s/r$ with $r, s \in \Z$ {satisfying} $0\leq |s| <r$ and $\gcd (r,s)=1$.
The case $\ta=\pi/2$ refers to the usual congruent numbers. We call such a triple $(a,b,c)$ a
\textit{rational $\ta$-triangle}.
A positive integer $N$ is called a \textit{$\ta$-congruent number} if there exists a rational $\ta$-triangle  with area  $N \sqrt{r^2-s^2}$, where $\ta$  is the angle between the sides $a$ and $b$; equivalently,  if there is a triple $(a, b, c)$ of (positive) rational numbers  satisfying
\begin{equation}
\label{eq2}
a^2+b^2-\frac{2s}{r}ab=c^2, \quad ab=2rN.
\end{equation}
A triple $(a,b,c)$ of rational numbers satisfying \eqref{eq2} it is called a \textit{rational $\ta$-triple for $N$}. It is also said that $(a,b,c)$ is a \textit{rational $\ta$-triangle for $N$},
if $a, b, c$ are all positive rational numbers. Note that if we have a  $\ta$-triple
$(a, b, c)$ for $N$, then $(|a|,|b|,|c|)$ is a rational $\ta$-triangle for $N$.
We identify the  $\ta$-triples $(a,b,c)$ and $(b,a,c)$ since the equations~\eqref{eq2} are symmetric in $a$ and $b$. It is clear that if a positive integer $N$ is a $\ta$-congruent number with a  $\ta$-triple
$(a ,b, c)$, then $N m^2$ is also a $\ta$-congruent number with the rational $\ta$-triple $(ma, mb, mc)$. Hence, one may  assume  the square-free positive numbers in the study of $\ta$-congruent numbers.

There are some natural questions related to the rational $\ta$-triples and hence
$\ta$-triangles for $\ta$-congruent numbers.
\begin{ques}
	\label{mq}
	Let $\ta\in (0,\pi)$ be an angle  with rational cosine and assume that  $N$ is a $\ta$-congruent number.
	\begin{itemize}
		\item [(i)] Is it possible to introduce new rational 	$\ta$-triples for $N$ using a given one? If yes, can one repeat this procedure infinitely many times?
		\item [(ii)]  Is it possible to find new rational $\ta$-triples for $N$ using given  two distinct ones?
	\end{itemize}
\end{ques}
The part (i) of Question~\ref{mq}, in the case $\ta=\pi/2$, was answered around $360$ years ago by Fermat \cite{Fermat,Hun}. Indeed, without any proof,  he provided an algorithm which is capable to produce infinitely many different right triangles for a given congruent number $N$ from a given right triangle.  In \cite{Halb-Hun}, Halbeisen and Hungerb\"{u}hler provided a detailed  proof for the Fermat's algorithm. They used their results to give an elementary proof of the fact that if $N$ is a congruent number and $P_0=(x_0, y_0)$ with $y_0\not =0$  is a rational point on the
\textit{congruent number elliptic curve} $E_N: y^2=x^3-N^2 x$, then $P_0$ cannot be of finite order;
in other words, it is a point of infinite order in the Mordell-Weil group $E_N(\Q)$ and so it does not belong to the set $T_N(\Q)=\left\{\infty, (0,0), (-N,0), (N,0)\right\}$, called the \textit{torsion subgroup of $E_N(\Q)$}, where $\infty=(0:1:0)$ denotes the \textit{point at infinity}.
In the recent work \cite{chan}, for $\ta=\pi/2$,  Chan gave a positive answer to the part (ii) of  Question~\ref{mq} using  geometric methods.

In \cite{fujw1,fujw2}, Fujiwara proved a theorem (see Theorem~\ref{fuj1}) on the \textit{$\ta$-congruent number elliptic curves} defined by
\begin{equation*}
\ent:~ y^2=x(x+(r+s)N)(x-(r-s)N),
\end{equation*}
where $r$ and $s$ are as mentioned earlier, {relating the $\ta$-congruent number problem to the existence of a non-2-torsion point in  $E_N^\ta(\Q)$.
	It is clear that the points $(0,0), (-(r+s)N,0)$ and $((r-s)N,0)$ are the only  $2$-torsion
	points in $E_N^\ta(\Q)$. This implies that $\Z/2\Z \times \Z/2\Z $ is contained in  $T_N^\ta(\Q)$, \textit{the torsion subgroup of $E_N^\ta(\Q)$}. Thus, based on the  Mazur's theorem~\ref{maz}
	on the classification of torsion subgroups of elliptic curves over $\Q$,
	the only possibilities for $T_N^\ta(\Q)$ are $\Z/2\Z \times \Z/2n\Z$ with $1\leq n \leq 4$.
	In \cite{fujw2}, Fujiwara determined certain conditions on  $\ta$ and $N$ so that $T_N^\ta(\Q)$
	is isomorphic to one of the above possibilities, see Theorem~\ref{fuj2}.
	We note that Fujiwara did not provide explicitly the points  of  $T_N^\ta(\Q)$ in his proof.}

In this  paper,  we give a positive answer to the both parts of Question~\ref{mq} for any  $\ta\in (0,\pi)$ with rational cosine and any $\ta$-congruent number $N$.
To do this,  in Section~\ref{ell-cong},
we  briefly review the arithmetic of $\ta$-congruent numbers and the Fujiwara's results.
Then,  in Section~\ref{FTALG}, we answer  to the part (i) of Question~\ref{mq} by generalizing  the Fermat's algorithm for the rational  $\ta$-triangles. The  Section~\ref{Add-dup} is devoted to show the
analogy between our generalization of the Fermat's algorithm and the duplication formula in  $ E_n^\ta(\Q)$.
We also show how to use the  addition  of two distinct points in $E_N^\ta(\Q)$ to give  a positive answer to  the part (ii) of Question~\ref{mq}.
In the last section, we provide  an alternative proof for the Fujiwara's theorem~\ref{fuj1}
and the right-to-left side of Theorem~\ref{fuj2}  using the  rational $\ta$-triples.
We also list all possible torsion points in $E_N^\ta(\Q)$ with corresponding rational $\ta$-triples.

\section{Elliptic curves and  $\ta$-congruent numbers}
\label{ell-cong}
An \textit{elliptic curve} over $\Q$ is a genus one smooth algebraic (plane) curve defined by the affine short Weierstrass form
$$E:~ y^2=x^3+ a x +b, \quad a,b \in \Z, \ 4 a^3+27 b^2 \not =0,$$
together with the point $\infty$.
The  set of rational points on $E$, denoted by $E(\Q)$, forms an Abelian group with respect to the chord-tangent  addition law.
To see more  on the geometry of elliptic curves,  we cite the reader to \cite{slm1}.
By  a celebrated  theorem of  Mordell and Weil \cite{slm1}, {the set $E(\Q)$ is} a finitely generated Abelian group with the identity $\infty$; in other words, we have
$$E(\Q)\cong  T(E,\Q) \times \Z^{r(E)},$$
where $T(E,\Q) $ is a finite group called the \textit{torsion subgroup} of $E(\Q)$ and
$r(E)$ is a non-negative integer called the \textit{Mordell-Weil rank} of $E(\Q)$.
Determining the Mordell-Weil rank of elliptic curves over $\Q$ is a very challenging problem since there is no  algorithm to determine it in general. But, in contrast,
Mazur  classified all possibilities for torsion subgroups of elliptic curves over $\Q$  by  studying  modular curves in a groundbreaking work \cite{mazur}.
Here is the resume of Mazur's results.
\begin{thm}[{Mazur}]
	\label{maz}
	Let $E$ be an elliptic curve over $\Q$.
	Then, the torsion subgroup $T(E,\Q)$ is one of the following finite groups:
	\begin{itemize}
		\item [\rm(i)]  $\Z/n\Z \ \text{with} \ 1 \leq n \leq 10 \  or \ n=12,$
		\item	[\rm(ii)]   $\Z/2\Z \times \Z/2n\Z \ \text{with} \ 1 \leq n \leq 4.$
	\end{itemize}
\end{thm}

Let us  recall the algebraic version of the addition law and the duplication
formula in $E_N^\ta(\Q)$.
The sum $P_0+P_1=(x_2,y_2)$, of  any two points $P_0=(x_0, y_0)$ and $P_1=(x_1, y_1)$ in  $\ent(\Q)$,
is given by the following rules:
\begin{itemize}
	\item[\rm(1)] If $x_0\neq x_1$, then $x_2={\lambda}^2-2sN-x_0-x_1$ and $ y_2=-(\lambda x_2 + \nu)$, where
	$$~\lambda =\frac{y_1-y_0}{x_1-x_0},~\text{and}~\nu =y_0 - \lambda x_0,$$
	\item[\rm(2)] If $P_0=P_1$, then
	$x_2={\lambda}^2-2sN-2x_0$ and $ y_2=-{\lambda}^3+(3x_0+2sN){\lambda}-y_0,$
	where
	$${\lambda}=\frac{3x_0^2+4sNx_0-(r^2-s^2)N^2}{2y_0}.$$
\end{itemize}
For the $\ta$-congruent number elliptic curve {$E_N^\ta$}, we denote by
$T_N^\ta(\Q)$  the torsion subgroup and
by {$r_N^{\ta}(\Q)$} the Mordell-Weil rank of $\ent(\Q)$.
In the recent decades, the notion of $\theta$-congruent numbers has attracted the interests of some  mathematicians, consult \cite{djps,fujw1,fujw2,j-s,kan1}. Among other results, Fujiwara \cite{fujw1} showed the following relation between $\theta$-congruent numbers and  {$r_N^{\ta}(\Q)$}.
\begin{thm}
	[Fujiwara]
	\label{fuj1}
	Let $\ta\in(0,\pi)$ be an angle with rational  cosine.
	\begin{itemize}
		\item[\rm(1)] A positive integer $N$ is a $\ta$-congruent number if and only if $\ent(\Q)$ has a point of order greater than $2$,
		\item[\rm(2)]  If $N\nmid 6$, then $N$ is a $\ta$-congruent number if and only if
		$r_N^{\ta}(\Q)>0$ for the group $\ent(\Q)$.
	\end{itemize}
\end{thm}

Moreover, in \cite{fujw2}, the torsion subgroups of $\theta$-congruent elliptic curves are classified by Fujiwara as follows.
\begin{thm}
	[Fujiwara]
	\label{fuj2}
	For $T_N^{\ta}(\Q)$, one of the following cases happens.
	\begin{itemize}
		\item[\rm{(1)}] $T_N^{\ta}(\Q)\cong \Z/2 \Z \times \Z/8\Z$ if and only if there exist integers $u, v>0$ such that $\gcd(u,v)=1$, $u$ and $v$ have opposite parity and satisfy either of the following:
		\begin{itemize}
			\item[\rm{(i)}] $N=1$, $r=8u^4v^4$, $r-s={(u^2-v^2)^4}$, $(1+\sqrt{2})v>u>v$,
			\item[\rm{(ii)}] $N=2$, $r=(u^2-v^2)^4$, $r-s=32u^4v^4$, $u>(1+\sqrt{2})v$,
		\end{itemize}
		\item[\rm{(2)}] $T_N^{\ta}(\Q)\cong \Z/2\Z\times \Z/6\Z$ if and only if there exist integers $u,v>0$ such that $\gcd(u,v)=1$, $u>2v$ and satisfy either of the following:
		\begin{itemize}
			\item[\rm{(i)}] $N=1$, $r=(u-v)^3(u+v)/2$, $r+s=u^3(u-2v)$,
			\item[\rm{(ii)}] $N=2$, $r=(u-v)^3(u+v)$, $r+s=2u^3(u-2v)$,
			\item[\rm{(iii)}] $N=3$, $r=(u-v)^3(u+v)/6$, $r+s=u^3(u-2v)/3$,
			\item[\rm{(iv)}] $N=6$, $r=(u-v)^3(u+v)/3$, $r+s=2u^3(u-2v)/3$,
		\end{itemize}
		\item[\rm{(3)}] $T_N^{\ta}(\Q)\cong \Z/2\Z\times \Z/4\Z$ if and only if
		\begin{itemize}
			\item[\rm{(i)}] $N=1$, $2r$ and $r-s$ are squares but not satisfy (1)(i),
			\item[\rm{(ii)}] $N=2$, $r$ and $2(r-s)$ are squares but not satisfy (1)(ii),
		\end{itemize}
		\item[\rm{(4)}] Otherwise, $T_N^{\ta}(\Q)\cong \Z/2\Z\times \Z/2\Z$.
	\end{itemize}
\end{thm}

Here, we recall the following one-to-one correspondence  between two sets with different objects which has crucial role in the proof of the above theorems and the results of the current paper.
For $\ta\in (0,\pi)$  with rational cosine and  a $\ta$-congruent number $N$, we consider the sets:
\begin{align*}
\Tc_N^\ta &= \left\{(a, b, c) \in \Q^3\setminus {(0,0,0)}:   c^2=a^2+b^2 -\frac{2s}{r}ab, ab=2rN \right\},\\
\Pc_N^\ta &=\left\{(x, y) \in \Q^2:  (x, y) \in   \ent(\Q)\setminus \{\infty\},  y \neq  0 \right\}.\end{align*}
Then, by a straightforward computation,  one  can easily check that
the sets $\Tc_N^\ta$ and $\Pc_N^\ta$  are bijective via the following maps:
\begin{align}
\phi :  (x, y) &\mapsto  \left(  \frac{y}{x}, \frac{2rxN}{y},
\frac{x^2 + (r^2-s^2) N^2}{y}\right), \label{tran2} \\
\psi :  (a,b,c) &\mapsto  \left(\frac{r N (a+c- sb/r)}{b},\frac{2 r^2 N^2 (a+c- sb/r)}{b^2} \right).\label{tran1}
\end{align}
Using the fact $a/2= r N/b$,  we can rewrite
$$\psi(a,b,c)=  \left(\frac{a(a+c- sb/r)}{2},\frac{a^2 (a+c- sb/r)}{2} \right).$$
Note that in the set  $\Pc_N^\ta$ it is  allowed  $a, b, c$ to be negative rational numbers and that $a, b$ have the same sign since their product is a positive integer.
Hence, all the four $\ta$-triples
$$(a,  b, c), ( a, b, - c), (- a,- b,  c), (- a,- b, -c)  \in \Pc_N^\ta $$
correspond to the rational $\ta$-triangle  $(|a|, |b|, |c|)$
for the $\ta$-congruent number $N$.
It is an easy fact that $a+c-sb/r >0$  for a given  rational $\ta$-triangle  $(a, b, c)$ for $N$. This means that any rational $\ta$-triangle  $(a, b, c)$ is mapped by $\psi$ on a point $(x, y)$ with positive coordinates.

We end this section by the following remark on $\ta$-triples.
\begin{rem}
	\label{crite}
	Using the equation~\eqref{eq2} it is easy to see that there exists a rational $\ta$-triple $(a, b, c)$
	with $|a|=|b|$  for  a $\ta$-congruent number $N$  if and only if the rational numbers 	
	$ 2 r(r - s) $  and  $ 2 r(r + s) $    are squares.
	Indeed,  we have
	$$
	\begin{array}{rl}
	a=b\Longleftrightarrow 2rN=a^2, & 2r(r-s)=(rc/a)^2, \vspace{.2cm} \\
	a=-b \Longleftrightarrow 2rN =a^2, & 2r(r+s)=(rc/a)^2.
	\end{array}
	$$
	%
	In this case, the  $\ta$-triple is given by
	$$(a, \mp a ,c)=(\sqrt{2 r N}, \sqrt{2 r N}, 2 \sqrt{(r \pm s) N}).$$
	In particular,  there exists a rational $\ta$-triangle  $(a, b, c)$ with $a=b=c$ for  a $\ta$-congruent number $N$  if and only if   $\ta=\pi/3$, $a=2$  and $N=1$.
\end{rem}

\section{Generalizing the Fermat's algorithm for rational $\ta$-triples}
\label{FTALG}
Let us recall the Fermat's algorithm \cite{Fermat,Hun} which gives infinitely many distinct right triangles for the congruent  numbers if only one of them exists.
\begin{thm}
	\label{fer0}
	Suppose that $N$ is a congruent number and  $(a_0, b_0, c_0)$ is
	a  rational right triangle for $N$.
	Then,
	\begin{align*}
	\label{eq3a}
	a_1&= \frac{2c_0 a_0 b_0}{b_0^2-a_0^2} = \frac{4 c_0N}{\sqrt{c_0^4 -16 N^2}},  \\
	b_1&= \frac{c_0^4 - 4 a_0^2 b_0^2}{2 c_0 (b_0^2-a_0^2)} = \frac{\sqrt{c_0^4 -16 N^2}}{2 c_0},\\
	c_1&= \frac{c_0^4 + 4a_0^2 b_0^2 }{2 c_0 (b_0^2-a_0^2)}= \frac{c_0^4 +16 N^2}{2 c_0 \sqrt{c_0^4 -16 N^2}},
	\end{align*}
	give  another rational  right triangle for $N$. Furthermore,  this procedure leads to infinitely many different	rational right triangles for $N$.
\end{thm}

First,  we adapt the Fermat's algorithm to the rational $\ta$-triples associated
to $\ta$-congruent numbers for any $\ta \in (0,\pi)$. The following theorem gives a naive algorithm for constructing new rational $\ta$-triples for the $\ta$-congruent numbers from given ones.
\begin{thm}
	\label{fer1}
	Suppose that $\ta \in (0,\pi)$ has a rational cosine and  $(a_0, b_0, c_0)$ is  a rational $\ta$-triple for a $\ta$-congruent number $N$. {If  $b_0^2- a_0^2 \neq 0 $, then }
	\begin{equation}
	\label{eq3}
	\begin{array}{l}
	a_1=\displaystyle
	{\frac{c_0^4- 4(r^2-s^2)a_0^2 b_0^2/r^2+ 4 s a_0 b_0  {c_0^2}/r}{2 c_0 (b_0^2-a_0^2)}}, \vspace{.2cm}\\
	b_1=\displaystyle \frac{4c_0^2 a_0 b_0}{2 c_0 (b_0^2-a_0^2)},  \vspace{.2cm}\\
	c_1=\displaystyle \frac{c_0^4 + 4(r^2-s^2)a_0^2 b_0^2/r^2 }{2 c_0 (b_0^2-a_0^2)},
	\end{array}
	\end{equation}
	is another rational $\ta$-triple  for $N$. It is different from $(a_0, b_0, c_0)$ provided that  $b_0^2- a_0^2 \neq  2 a_0 c_0$.
	In particular, $(|a_1|, |b_1|, |c_1|)$ is a  rational $\ta$-triangle for $N$.
\end{thm}
\proof
Given any $m, n \in \Q^*$, defining $X, Y $ and $Z$ as
$$X=m^2-\left(\frac{r^2-s^2}{r^2}\right)n^2 + \frac{2smn}{r}, \quad Y =2 m n, \quad Z= m^2+\left(\frac{r^2-s^2}{r^2}\right)n^2,$$
one can get a  rational $\ta$-triple for $N_0:=|{XY}/{(2r)}|$.
It is an easy task to check out the relation $X^2+Y^2- 2sXY/r=Z^2$.
Setting $m=c_0^2$ and $n= 2 a_0 b_0$  implies that
\begin{align*}
X&= c_0^4- 4\left(\frac{r^2-s^2}{r^2}\right) a_0^2 b_0^2 + 4\left(\frac{s}{r}\right)   {c_0^2} a_0 b_0, \\ Y&=4c_0^2 a_0 b_0, \\
Z&= c_0^4 + 4\left(\frac{r^2-s^2}{r^2}\right)  a_0^2 b_0^2.
\end{align*}
Thus,  we have
$$N_0=\left|\frac{XY}{2r}\right|= \left |\frac{2 c_0^2 a_0 b_0}{r} \left(c_0^2 - \frac{2(r-s)a_0b_0}{r} \right)
\left(c_0^2 + \frac{2(r+s)a_0b_0}{r} \right) \right |.$$
Since $c_0^2=a_0^2+b_0^2 -2sa_0b_0/r$, it follows
$N_0=|2  a_0 b_0 c_0^2 (a_0^2-b_0^2)^2/r|$,  which is nonzero by the assumption.
Now, we define  $(a_1, b_1, c_1):=(X/d_0, Y/d_0, Z/d_0),$
where $d_0=2c_0 (b_0^2-a_0^2)\neq 0$, again by the same assumption.
Clearly, we have  $a_1^2+b_1^2 -2sa_1b_1/r=c_1^2$ and
{$$ \frac{a_1 b_1}{2 r}=\frac{X Y}{2r  4 c_0^2(b_0^2-a_0^2)^2}
	=\frac{2 r N_0}{2 r  4 c_0^2\left( b_0^2-a_0^2\right) ^2}
	= \frac{4  a_0 b_0  c_0^2 \left( b_0^2-a_0^2\right)^2}{2 r 4 c_0^2
		\left( b_0^2-a_0^2\right)^2}=\frac{a_0 b_0}{2 r}.$$}
Thus,  both the $\ta$-triangles   $(|a_0|, |b_0|, |c_0|)$ and
$(|a_1|, |b_1|, |c_1|)$  have the same area $N\sqrt{r^2-s^2}$. We note that if $2 a_0 c_0 = b_0^2 -a_0^2$, then  $b_1=b_0$ and using
$c_0^2=a_0^2+ b_0^2 - 2s a_0 b_0/r$  one can get
$3 a_0^4-8 s a_0^3 b_0/r + 6 a_0^2 b_0^2 - b_0^4 =0$ which
divides both of the equations $a_1-a_0=0$ and $c_1-c_0=0$.
Thus,  by the assumption we have $(a_1, b_1, c_1)\neq  (a_0, b_0, c_0)$.
Therefore, the proof of Theorem~\ref{fer1} is completed.

Note that in  the above  theorem it is possible to have $a_1=b_1$.
Indeed, substituting
$c_0^2= a_0^2+ b_0^2- 2 s a_0 b_0 /r$ in the equations~\eqref{eq3} we have  $a_1=b_1$ if and only if 	$(a_0-b_0)^4=32 r (r-s)N^2$.
We refer the reader to  Table~\ref{tab:1} in the end of  this paper to see examples satisfying these conditions.

The following theorem shows that, under  mild conditions, repeating the procedure given by the theorem~\ref{fer1}
leads to  infinitely many distinct rational $\ta$-triangles for a given $\ta$-congruent number.
These conditions will be checked and used in the second half of Section 5 when Fujiwara's Theorem~\ref{fuj1} is proved.

\begin{thm}
	\label{fer2}
	For any  integer $n\geq 0$, suppose that
	$(a_{n+1}, b_{n+1}, c_{n+1})$  is the rational $\ta$-triple for  $N$
	obtained by  Theorem~\ref{fer1} from $(a_n, b_n, c_n)$   satisfying
	the conditions  $b_n^2-a_n^2\neq 0 $ and  $b_n^2- a_n^2 \neq  2 a_n c_n$.
	Then,  $|c_n|\not =|c_{n'}|$ for any two distinct non-negative integers $n$, $n'$.
\end{thm}
\proof
Fix an arbitrary  number $n\geq 0$. Since  $N=\left|\frac{a_n b_n}{2r}\right|$, we have
\begin{equation}
\label{eq4}
a_n^2 b_n^2=4 r^2 N^2.
\end{equation}
On the other hand, the equality $a_n^2+b_n^2- 2s a_n b_n/r=c_n^2$ implies that
$$c_n^4- 8 r N \left(a_n- \frac{s b_n}{r}\right) \left(b_n- \frac{s a_n}{r}\right) =\left(b_n^2-a_n^2\right)^2> 0.$$
Therefore,
\begin{equation}
\label{eq4a}
d_n:=b_n^2-a_n^2=
\sqrt{c_n^4- 8 r N \left(a_n- \frac{s b_n}{r}\right)\left(b_n- \frac{s a_n}{r}\right)}.
\end{equation}
Using the equations~\eqref{eq3}, \eqref{eq4}, and doing  a little bit of algebraic simplifications one can get that
\begin{equation}
\label{eq5}
d_n=\sqrt{c_n^4 + 8 s Nc_n^2- 16 N^2(r^2 - s^2)}, \quad
|c_{n+1}|=\frac{c_n^4+ 16 N^2(r^2-s^2)}{2 d_n |c_n|}.
\end{equation}
Note that $d_n$ is a rational number because  $c_{n+1} \in \Q^*$. Assume that $|c_n|=u/v$, where $u, v$ are in lowest term. To prove the assertion of the theorem, we consider the following two cases.\\
{\bf Case 1:  $u$ is even.}
Let us write  $u= 2^t\cdot \tilde{u}$, where $t \geq 1$ and $\tilde{u}$ is an odd number, hence
$c_n=2^t\cdot \tilde{u}/v$ and $N=2^{\ell} \cdot \tilde{N},$ where  $\ell \geq 1$ and
$\tilde{N}$ is an odd number.
Thus, we have
\begin{equation}
\label{8}
\begin{array}{l}
\displaystyle d_n^2 =\frac{2^{4t}\cdot\tilde{u}^4 + 2^{2t+{\ell}+3}\cdot\tilde{N}sv^2\tilde{u}^2- 2^{2{\ell}+4} \cdot \tilde{N}^2 (r^2-s^2)v^4}{v^4}=\left(\frac{2^{m} u_0}{v^2}\right)^2,\vspace{.2cm}\\
\displaystyle c_n^4 + 16 N^2(r^2 - s^2)=\frac{2^{4t}\cdot\tilde{u}^4+2^{2{\ell}+4}\cdot\tilde{N}^2(r^2-s^2)v^4}{v^4}=\left(\frac{2^{m} u_1}{v^2}\right) ^2,
\end{array}
\end{equation}
where the both numbers $u_0, u_1$ are odd and
$2 \leq m \leq 2t$. Therefore, from the equations~\eqref{eq5} and \eqref{8}, we can write
$$|c_{n+1}|= \frac{2^{m-t-1}\cdot u'}{v'},$$
where $u_0, u_1, u', v'$	 are odd numbers.
Since $m< 2t+1$, we have $ m-t-1<t$, which implies that:
\textit{If $|c_n|={2^t\cdot \tilde{u}}/{v}$ with odd $\tilde{u}, v$ and $t\geq0$, then we have $|c_{n+1}|={2^{t'}\cdot u'}/{v'}$ with odd $u', v'$ and $0 \leq t' <t$.} \\
{\bf Case 2:  $u$ is odd.}
We write $v=2^t\cdot \tilde{v}$, where $t\geq 0$ and $\tilde{v}$ is odd, hence
$|c_n|=u/(2^t\cdot\tilde{v})$ from which we obtain
$$ d_n^2= \left(\frac{\tilde{u}}{2^{2t}\tilde{v}^2}\right)^2, \quad  c_n^4 + 16 N^2(r^2 - s^2) = \frac{\bar{u}}{2^{4t}\cdot \tilde{v}^4}$$
for some odd positive numbers $\bar{u}$ and $ \tilde{u}$.
After some simple computations, we have
$$|c_{n+1}|=\frac{u'}{2^{t+1}\cdot  v'},$$
where $ u', v'$ are odd integers and $\gcd(u',v')=1.$ Thus, we conclude that:

\textit{If $|c_n|={u}/{2^t\cdot \tilde{v}}$ with odd $u, \tilde{v}$, and $t\geq0$, then $|c_{n+1}|={u'}{2^{t+1}\cdot v'}$ with odd $u', v'$ and $0 \leq t' <t$.}

Therefore, we have finished the proof of Theorem~\ref{fer2}.

Considering the equations~\eqref{eq4a} and  \eqref{eq5} in  the proof of Theorem~\ref{fer2}, we can reformulate Theorem~\ref{fer1} as follows.
\begin{cor}
	\label{fer3}
	Suppose that $\ta \in (0, \pi)$ has a rational cosine and  $(a_0, b_0, c_0)$ is  a rational $\ta$-triple
	with  $d_0\neq 0 $ and $d_0\neq 2a_0c_0$ for a $\ta$-congruent number $N$. Then,
	\begin{equation*}
	(a_1, b_1, c_1) =\left(\frac{d_0}{2 c_0}, \frac{4 r c_0 N}{d_0},\frac{c_0^4 + 16N^2(r^2 - s^2)}{2 c_0 d_0}\right)
	\end{equation*}
	is  a distinct rational  $\ta$-triple for $N$.
\end{cor}

\section{Rational $\ta$-triples and the addition  in  $\ent(\Q)$}
\label{Add-dup}

In this section, we show that our Fermat-type algorithm,  given in Theorem~\ref{fer1}, is essentially  doubling of points  in $E_N^\ta(\Q)$.
\begin{thm}
	\label{dup}
	Let  $\ta \in (0,\pi)$ be an angle with  rational cosine.  Given  a rational $\ta$-triple
	$(a_0, b_0, c_0)$ with $ d_0 \neq 0 $ and $d_0\neq  2 a_0 c_0$ for  a $\ta$-congruent number $N$, we assume that
	$(a_1, b_1, c_1)$ is a $\ta$-triple for $N$ obtained  by Theorem~\ref{fer1}.
	Moreover, we   let $(x_i,y_i)$  be the rational points in $E_N^\ta(\Q)$ corresponding to
	$(a_i, b_i, c_i)$ {by the map~\eqref{tran2}} for  $i=0,1$. Then,
	$$[2](x_0, y_0)=(x_1, -y_1).$$
\end{thm}
\proof
By Corollary~\ref{fer3},  the rational $\ta$-triple $(a_1, b_1, c_1)$
is given by
$$	(a_1, b_1, c_1) =\left(\frac{d_0}{2 c_0}, \frac{4 r c_0 N}{d_0},
\frac{c_0^4 + 16N^2(r^2 - s^2)}{2 c_0 d_0}\right), $$
which is corresponded by \eqref{tran1} to the point $(x_1, y_1)= (c_0^2/4, c_0 d_0/8)$ in $E_N^\ta(\Q)$. Indeed,  replacing $N= a_0 b_0/2 r$ and
the above values for $a_1, b_1$ and $c_1$, and doing some simplifications lead to
\begin{align*}
x_1 &= \frac{r N (a_1+ c_1 - s b_1/r)}{b_1} = \frac{c_0^2}{4},\\
y_1 &= \frac{2 r^2 N^2 (a_1+ c_1 - s b_1/r)}{b_1^2} = \frac{c_0 d_0}{8}.
\end{align*}
On the other side, replacing $N= a_0 b_0/2 r$ and using \eqref{tran1},
the triple  $(a_0, b_0, c_0)$  corresponds to  $(x_0,y_0)$ in $E_N^\ta(\Q)$  with coordinates,
$$x_0=\frac{a_0 (a_0+ c_0 - s b_0/r)}{2},\quad
y_0= \frac{a_0^2 (a_0+ c_0 - s b_0/r)}{2}.$$
Using the last equations and  applying   the duplication formula described in Section~2  one can show that
for  $(x_2,y_2):=[2](x_0,y_0)$ we have,
\begin{align*}
x_2&={\lambda}^2-2Ns-2x_0 = {c_0^2}/{4}= x_1,\\
y_2&=-{\lambda}^3+(3x_0+2Ns){\lambda}-y_0=-\frac{c_0 d_0}{8}= - y_1,
\end{align*}
where $\lambda$ is given by
\begin{align*}
{\lambda}=\frac{3x_0^4+4sN x_0^2-4 N^2(r^2-s^2)}{2 y_0}= \frac{2 a_0+c_0}{2}.
\end{align*}
Thus, we have completed the proof of Theorem~\ref{dup}.

In the following,  using the addition law in  $E_N^\ta(\Q)$, we provide a way to find new rational $\ta$-triples for $N$   from  given two distinct  ones.

Let $(a_0, b_0, c_0)$ and $(a_1, b_1, c_1)$ be two distinct $\ta$-triples for a $\ta$-congruent number $N$, where $\ta \in (0,\pi )$ has rational cosine as before.
Then, applying the transformation~\eqref{tran1}, we obtain the following two points $P_i=(x_i,y_i)$ in $E_N^\ta(\Q)$ for  $i=0,1$, where
\begin{equation*}
x_i= \frac{a_i ( a_i+ c_i)}{2}- s N, \quad y_i= \frac{a_i^2 ( a_i+ c_i)}{2}- s N a_i.
\end{equation*}
Defining $t_i=a_i( a_i + c_i)$ for  $i=0,1$   and
\begin{equation*}
\lambda=\frac{y_1-y_0}{x_1-x_0} =\frac{a_0 t_0- a_1 t_1 - 2 s N (a_0-a_1)}{  t_0 -  t_1}
\end{equation*}
we obtain that the point $P_2=(x_2, y_2):=P_0+ P_1$ has the coordinates
\begin{equation*}
\begin{array}{l}
x_2=\displaystyle\frac{(a_0-a_1)^2( t_0- 2 sN)( t_1 - 2 s N)}{( t_0 - t_1)^2},
\vspace{.2cm}\\
y_2=\displaystyle\frac{(a_0-a_1) ( t_0- 2 sN)( t_1 - 2 s N) T}{(  t_0 -  t_1)^3},
\end{array}
\end{equation*}
where
\begin{equation}
\label{t5}
T=4 r^2 N^2 +  t_0 t_1+ a_0 a_1 \left( t_0+ t_1 - 4 s N\right).
\end{equation}

Therefore, using the transformation~\eqref{tran2} and some algebraic computations, we conclude the following  theorem.

\begin{thm}
	\label{fer6}
	Suppose that $\ta \in (0,\pi)$ has a rational cosine, $(a_0, b_0, c_0)$ and  $(a_1, b_1, c_1)$  are two distinct rational $\ta$-triples for a $\ta$-congruent number $N$. 	Then,
	the $\ta$-triple	$(a_2, b_2, c_2)$ given by
	\begin{equation*}
	\begin{array}{l}
	a_2  =\displaystyle \frac{4 r^2 N^2 +  t_0 t_1+ a_0 a_1 \left( t_0+ t_1 - 4 s N\right)}{(a_0-a_1) (  t_0 -  t_1)}, \quad b_2= \frac{2 r N}{a_2}, \vspace{.2cm}\\
	c_2  =\displaystyle \frac{(a_0-a_1)^4( t_0- 2 sN)^2( t_1 - 2 s N)^2 + (r^2-s^2)(t_0 - t_1)^2 N^2   }{(a_0-a_1) (t_0 - t_1) ( t_0- 2 sN)( t_1 - 2 s N) T},
	\end{array}
	\end{equation*}
	is another rational $\ta$-triple for $N$, where $T$ is given by \eqref{t5}.
	In particular, $(|a_2|, |b_2|, |c_2|)$ is a rational $\ta$-triangle for $N$.
\end{thm}

\proof

We may   assume that $a_0\neq a_1$ to have two distinct triples. Then,
the relations $a_i b_i= 2 r N$  for $i=0,1$  give us $a_0\neq b_1$ and hence $c_0\neq c_1$. Thus, we have $t_0 \neq t_1$, hence both  $x_2$  and $y_2$ are well-defined.
Moreover, none of the factors in the defining equations of $x_2$  and $y_2$ can be zero, since the distinctness of $\ta$-triples implies $x_i\neq 0$ and $t_i\neq 2 s N$ for $i=0,1$.

In order to show that $(a_2, b_2, c_2)$ is a rational $\ta$-triangle for $N$, it is enough to check that  $a_2 b_2 = 2 r N$ and $a_2^2+ b_2^2 -  2 s a_2 b_2/r = c_2^2$. The first equality is trivial and the second one is equivalent to $a_2^2+ 4 r^2 N^2/a_2^2 -  4 s N - c_2^2=0$ by considering the first one. It is straightforward to check the latter equality, so we leave it to the reader. Therefore, we have a new rational $\ta$-triangle for $N$.

\section{Proof of Fujiwara's theorems}
\label{fujproof}

First,  using the equations \eqref{eq2}, \eqref{tran1} and the duplication formula, we prove the right-to-left  side of Theorem~\ref{fuj2} case by case as  follows. For the other side, we refer the reader to Fujiwara's  original proof in \cite{fujw2}.

\begin{itemize}
	\item[] Case (1)(i). $N=1$, $r=8u^4v^4$, $r-s=(u^2-v^2)^4$, $(1+\sqrt{2})v>u>v>0$,
	{and} $\gcd(u,v)=1$.
\end{itemize}
From $ab=2rN$, letting  $a=4uv^3$ and $b=4u^3v$, we  get $c=2(u^4-v^4)$ satisfying $c^2=a^2+b^2-2abs/r$. By \eqref{tran1}, the  $\ta$-triple $(a,b,c)$ transforms to the point $P_1=(x_1, y_1)$ or $-P_1$ with
$$x_1=(u+v)(u^2+v^2)(-v^2+2uv+u^2)(u-v)^3,\quad y_1=4uv^3x_1,$$
which is of order~$8$, carried out by {\sf{SAGE}}. By the duplication formula, we get the other order~$8$ points $P_j=(x_j, y_j)$ and $-P_j$, $j=2,3,4$, where
\begin{align*}
x_2&= (u-v)(u^2+v^2)(-v^2-2uv+u^2)(u+v)^3, \quad y_2=4uv^3x_2, \\
x_3&= (u-v)(u^2+v^2)(-v^2+2uv+u^2)(u+v)^3, \quad y_3=4u^3vx_3,\\
x_4&= (u+v)(u^2+v^2)(-v^2-2uv+u^2)(u-v)^3, \quad y_4=4u^3vx_4,
\end{align*}
so that $[2](\pm P_k)=((u^4-v^4)^2,\pm 4u^2v^2(u^4-v^4)^2)$ and $[4](\pm P_k)=((u^2-v^2)^4,0)$ for $k=1,2,3,4$. Hence, $T_N^{\ta}(\Q)\cong \Z/2 \Z \times \Z/8\Z$.
\begin{itemize}
	\item[] Case (1)(ii). $N=2$, $r=(u^2-v^2)^4$, $r-s=32u^4v^4$, $u>(1+\sqrt{2})v>0$, {and} $\gcd(u,v)=1$.
\end{itemize}
The $\ta$-triple $(a,b,c)=(2(u+v)(u-v)^3,2(u+v)^3(u-v),8uv(u^2+v^2))$ satisfies \eqref{eq2} and hence, by \eqref{tran1}, transforms to the order~$8$ point $P_1=(x_1,y_1)$ or $-P_1$ with
$$x_1=16uv^3(u^2+v^2)(u^2+2uv-v^2),\quad y_1=2(u+v)(u-v)^3x_1.$$
By the duplication formula, we get the other order~$8$ points $P_j=(x_j, y_j)$ and $-P_j$, $j=2,3,4$, having the coordinates
\begin{align*}
x_2&=16u^3v(u^2+v^2)(v^2+2uv-u^2), \quad y_2=2(u+v)(u-v)^3x_2, \\
x_3&=16u^3v(u^2+v^2)(u^2+2uv-v^2),\quad  y_3=2(u-v)(u+v)^3x_3,\\
x_4&=16uv^3(u^2+v^2)(v^2+2uv-u^2),\quad  y_4=2(u-v)(u+v)^3x_4.
\end{align*}
It is easy to check that $[2](\pm P_k)= (16u^2v^2(u^2+v^2)^2,\pm 32u^2v^2(u^4-v^4)^2)$ and  $[4](\pm P_k)=(64u^4v^4,0)$ for $k=1,2,3,4$. Hence, $T_N^{\ta}(\Q)\cong \Z/2 \Z \times \Z/8\Z$.
\begin{itemize}
	\item[] Case (2)(i). $N=1$, $r=(u-v)^3(u+v)/2$, $r+s=u^3(u-2v)$, $u>2v>0$, $\gcd(u,v)=1$, and
\end{itemize}
\begin{itemize}
	\item[] Case (2)(iii). $N=3$, $r=(u-v)^3(u+v)/6$, $r+s=u^3(u-2v)/3$, $u>2v>0$, $\gcd(u,v)=1$.
\end{itemize}
Considering $(a,b,c)=(u^2-v^2,(u-v)^2,2uv)$ which satisfies \eqref{eq2}, we obtain, by \eqref{tran1}, the order~$6$ point $P_1=(x_1,y_1)$ or $-P_1$ with
$$x_1=(2u-v)u^2v, \quad y_1=(u^2-v^2)x_1.$$
By the duplication formula, we get the other order~$6$ points $P_j=(x_j,y_j)$ and $-P_j$ for $j=2,3$ with coordinates
$$
x_2=(2v-u)uv^2,\quad x_3=(2u-v)(2v-u)uv,
$$
and
$$
y_2=(u^2-v^2)x_2, \quad y_3=(u-v)^2x_3,
$$
so that $[2](\pm P_k)=(u^2v^2,\pm u^2v^2(u-v)^2)$ for $k=1,2,3$.
Hence, we have $T_N^{\ta}(\Q)\cong \Z/2 \Z \times \Z/6\Z$.
\begin{itemize}
	\item[] Case (2)(ii). $N=2$, $r=(u-v)^3(u+v)$, $r+s=2u^3(u-2v)$, $u>2v>0$, $\gcd(u,v)=1$, and
\end{itemize}
\begin{itemize}
	\item[] Case (2)(iv). $N=6$, $r=(u-v)^3(u+v)/3$, $r+s=2u^3(u-2v)/3$, $u>2v>0$, $\gcd(u,v)=1$.
\end{itemize}
The $\ta$-triple $(a,b,c)=(2(u^2-v^2),2(u-v)^2,4uv)$  satisfies \eqref{eq2} and, by \eqref{tran1},  transforms to the order~$6$ point $P_1=(x_1,y_1)$ or $-P_1$ with $x_1=4(2u-v)u^2v$, $y_1=2(u^2-v^2)x_1$.
By the duplication formula, we get the other order~$6$ points $P_j=(x_j,y_j)$ and $-P_j$, $j=2,3$, with coordinates
$$
x_2=4(2v-u)uv^2,\quad
x_3=4(2u-v)(2v-u)uv,
$$
and
$$
y_2=2(u^2-v^2)x_2,\quad
y_3=2(u-v)^2x_3,$$
so that $[2](\pm P_k)=(4u^2v^2,\pm 8u^2v^2(u-v)^2)$. Hence, $T_N^{\ta}(\Q)\cong \Z/2 \Z \times \Z/6\Z$.

\begin{itemize}
	\item [] Case 3(i). $N=1$, $2r$ and $r-s$ are squares but not satisfy (1)(i).
\end{itemize}
We may let $r=2u^2$ and $s=2u^2-v^2$ and get the $\ta$-triple $(a,b,c)=(2u,2u,2v)$ which satisfies \eqref{eq2} and, by \eqref{tran1}, transforms to the order~$4$ point
$P_1=((2u+v)v,2uv(2u+v))$ or $-P_1$. By the duplication formula, we get the other order~$4$ points
$P_2=((v-2u)v,2uv(v-2u))$ and $-P_2$,
so that $[2](\pm P_j)=(v^2,0)$ for $j=1,2$. Hence, $T_N^{\ta}(\Q)\cong \Z/2 \Z \times \Z/4\Z$.
\begin{itemize}
	\item [] Case (3)(ii). $N=2$, $r$ and $2(r-s)$ are squares but not satisfy (1)(ii).
\end{itemize}
We may let $r=u^2$ and $s=u^2-2v^2$ and get the $\ta$-triple $(a,b,c)=(2u,2u,4v)$ which satisfies \eqref{eq2} and transforms to the order~$4$ point $P_1=(4(u+v)v,8(u+v)uv)$  or $-P_1$.
By the duplication formula, we get the other order~$4$ points $P_2=(4(v-u)v,8(v-u)uv)$ and $-P_2$, so that $[2](\pm P_j)=(4v^2,0)$ for $j=1,2$. Hence, $T_N^{\ta}(\Q)\cong \Z/2 \Z \times \Z/4\Z$.
\begin{itemize}
	\item[] Case (4). Otherwise.
\end{itemize}
In this case, it is evident that $T_N^{\ta}(\Q)\cong \Z/2 \Z \times \Z/2\Z$ because the image of $\phi$ is the set of all rational points on $E^{\ta}_N(\Q)$, except for the three points on the $x$-axis, which are the non-trivial points of order~$2$, and for the point at infinity.

In Table~\ref{tab:1}, we bring in the possible corresponding rational $\ta$-triples to the above torsion points of orders~$4$, $6$, $8$, obtained by the transformation~\eqref{tran2} and Theorem~\ref{fer1}.
Therefore, we have completed the proof of the right-to-left side of Fujiwara's Theorem~\ref{fuj2}.

Secondly,  we note that
given any $\ta$-congruent  number $N$ and  rational $\ta$-triangle $(a, b,c)$ one has
$a+c -s b/r > 0$ which implies that  its image $\psi (a,b,c)$,  by the transformation~\eqref{tran1},  has a nonzero   $y$-coordinate.
Hence, it is a non-2-torsion point in $E_N^\ta(\Q)$  as desired.
Conversely, by the fact that  the only $2$-torsion points in $E_N^\ta(\Q)$ are $(0,0)$, $(-(r+s) N, 0)$,  $((r-s) N, 0)$ and the point at the infinity, any non-$2$-torsion
gives us  a rational $\ta$-triangle for $N$ by transformation~\eqref{tran2}.
This shows the part (1) of the  Fujiwara's theorem~\ref{fuj1}.

Finally, if $N\nmid 6$ is a $\ta$-congruent number,  $(x_0, y_0) \in E_N^{\ta}(\Q)$
with $y \not = 0$ and
$$(a_0, b_0, c_0)=\phi(x_0, y_0)=\left(\frac{y_0}{x_0},  \frac{2r x_0 N}{y_0},
\frac{x_0^2 + (r^2-s^2) N^2}{y_0}\right)$$
is its corresponding rational $\ta$-triangle, then it
satisfies  the conditions of Theorems~\ref{fer1} and \ref{fer2}.
Indeed,  we have  $b_0^2= a_0^2$ if and only if
$a_0=  \pm b_0$ if and only if $y_0^2= \pm 2 r x_0^2 N$.  On substituting
$y_0^2= x_0 (x_0^2 + 2 s N x_0 - (r^2-s^2)N^2)$ and some simplification we obtain that
$a_0=  \pm b_0$ if and only if $$x_0^2 - 2 N (r \pm s) x_0- (r^2-s^2) N^2=0,$$ which is soluble in $x_0 \in \Q$  if and only if $ 2 r (r-s)$ is a rational square.
But, the proof of the left-to-right side of Fujiwara's Theorem~\ref{fuj2} in \cite{fujw2} shows the condition holds only in the cases $N= 1,2,3$ and $6$, each of which gives us a contradiction and hence
we have $a_0^2 \neq b_0^2$.

On the other hand, a simple algebraic calculation shows that for the above $\ta$-triple
$(a_0, b_0, c_0)$, the relation $b_0^2- a_0^2= 2 a_0 c_0$ holds
if and only if the following quartic equation
$$3\,{x_0}^{4}+8\,Ns{x_0}^{3}-  6\,{N}^{2} \left( {r}^{2}-{s}^{2}
\right) {x_0}^{2}-{N}^{4}({r}^{2}-{s}^{2})^2=0$$
has  rational solution $x_0$. This is equivalent to saying that the $x$-coordinate of $[2]P_0$ is equal to $x_0$; in other words,
$P_0$ is a point of order~$3$. But,  this is impossible by the Fujiwara's theorem~\ref{fuj2} since $N \nmid 6$. Thus, we have  $b_0^2- a_0^2\neq 2 a_0 c_0$ as desired.

Now, let  $P_0=(x_0, y_0)$ and for any $n\geq 0$ suppose that $(a_{n+1}, b_{n+1}, c_{n+1})$ is the $\ta$-triple obtained  by applying Theorem~\ref{fer1} for $(a_{n}, b_{n}, c_{n})$ which are all distinct. Denote by  $P_n=(x_n,y_n)$ the rational point in $E_N^\ta(\Q)$ corresponding  to  $(a_{n}, b_{n}, c_{n})$. Then,  a similar argument as given in  the proof of Theorem~\ref{dup}
shows that $x_n=c_n^2/4$  and  $P_n=[2^n] P_0$  for any $n\geq 0$.
Since $|c_n|\neq |c_{n'}|$, the points  $P_n $, $ P_{n'}$ are distinct  in $E_N^\ta(\Q)$ for $n\neq n' \geq 0$ by Theorem~\ref{fer2}.
In other words, the point $P_0$ is of infinite order so the rank of $E_N^\ta(\Q)$ is greater than  or equal to one. The converse is trivial  by considering the Fujiwara's theorem~\ref{fuj2} on the possible torsion points in $E_N^\ta(\Q)$.
Therefore, we have completed the proof of Theorem~\ref{fuj1}.

\begin{table}[h]
	\caption{Rational $\ta$-triples obtained by \eqref{tran2} and Theorem~\ref{fer1} for the  points of $T_N^\ta(\Q)$}
	\label{tab:1}
	{\small	\begin{tabular}{cllll}
			\hline\noalign{\smallskip}
			Cases &  Points  & $(a_0,b_0,c_0)$  &  $(a_1,b_1,c_1)$ & $\cdots$ \\
			\noalign{\smallskip}\hline\noalign{\smallskip}
			{(1)(i)}  & $\pm P_1, \pm P_2$ & $(4uv^3,4u^3v,2(u^4-v^4))$ & $(4u^2v^2,4u^2v^2,2(u^2-v^2)^2)$ & {*}\\
			& $\pm P_3, \pm P_4$ & $(4u^3v,4uv^3,2(u^4-v^4))$ &  $(4u^2v^2,4u^2v^2,2(u^2-v^2)^2)$ & {*} \\
			& $[2](\pm P_k)$ &  $(4u^2v^2,4u^2v^2,2(u^2-v^2)^2)$  & & {*} \\
			\noalign{\smallskip}\hline\noalign{\smallskip}
			{(1)(ii)} & $\pm P_1, \pm P_2$ & $a_0=2(u+v)(u-v)^3$ &  $a_1=2(u^2-v^2)^2$  & \\
			&  & $b_0=2(u+v)^3(u-v)$ & $b_1=2(u^2-v^2)^2$ & {*}\\
			& & $c_0=8uv(u^2+v^2)$  & $c_1=16u^2v^2$& \\
			& $\pm P_3, \pm P_4$ & $a_0=2(u-v)(u+v)^3$ & $a_1=2(u^2-v^2)^2$ &\\
			&  & $b_0=2(u-v)^3(u+v)$ & $b_1=2(u^2-v^2)^2$ & {*}\\
			&  & $c_0=8uv(u^2+v^2)$ & $c_1=16u^2v^2$ &\\
			& $[2](\pm P_k)$ & $a_0=2(u^2-v^2)^2$ & & \\
			&                        & $b_0=2(u^2-v^2)^2$  & & {*}\\
			&                        & $c_0=16u^2v^2$ &  & \\
			\noalign{\smallskip}\hline\noalign{\smallskip}
			{(2)(i),} & $\pm P_1, \pm P_2$ & $(u^2-v^2,(u-v)^2,2uv)$ & $(u^2-v^2,(u-v)^2,2uv)$ & {**} \\
			{(2)(iii)} & $\pm P_3, [2](\pm P_k)$ & $((u-v)^2,u^2-v^2,2uv)$ &  $((u-v)^2,u^2-v^2,2uv)$ & {**} \\
			\noalign{\smallskip}\hline\noalign{\smallskip}
			{(2)(ii),} & $\pm P_1, \pm P_2$ & $(2(u^2-v^2),2(u-v)^2,4uv)$ &  $(2(u^2-v^2),2(u-v)^2,4uv)$ &{**} \\
			{(2)(iv)} & $\pm P_3, [2](\pm P_k)$ & $(2(u-v)^2,2(u^2-v^2),4uv)$ &  $(2(u-v)^2,2(u^2-v^2),4uv)$ & {**} \\
			\noalign{\smallskip}\hline\noalign{\smallskip}
			{(3)(i)}  & $\pm P_1, \pm P_2$ & $(2u,2u,2v)$  & & {*}\\
			\noalign{\smallskip}\hline\noalign{\smallskip}
			{(3)(ii)} & $\pm P_1, \pm P_2$ & $(2u,2u,4v)$  & & {*} \\
			\noalign{\smallskip}\hline
	\end{tabular}}
\end{table}

{The symbol ``$*$" in the above table shows that the theorem~\ref{fer1} cannot be  applied further because of the lack of the condition $b_i^2-a_i^2\neq 0$ for some index $i=1$ or $2$ depending on the cases. By the symbol ``$**$", we mean that applying Theorem~\ref{fer1} leads to a
	rational 	$\ta$-triple which can be identified with the given one if we permute  the $a$ and $b$ components.
	This happens  in the cases $\rm (2)(i)$ to $\rm (2)(iv)$, since the condition $b_0^2-a_0^2= 2 a_0 c_0$ does not hold.}



\begin{thebibliography}{99}



\bibitem {chan} \textsc{S. Chan,} \textit{Rational right triangles of a given area},
The American Mathematical Monthly, \textbf{125}(8) pp. 689--703.


\bibitem {djps}
A. Dujella,  A. S. Janfada,  C. J. Peral,  and  S. Salami,
\textit{On the  high rank $\pi/3$ and $2\pi/3$-congruent number elliptic curves,}
Rocky Mountain J.  Math. \textbf{44}(6) (2014) 1867--1880.

\bibitem{Fermat}\textsc{Fermat, P.}:\ \emph{Fermat's Diophanti Alex. Arith., 1670 in Ouvres III},
(Minist\`{e}re de l'instruction publique, ed.),
Gauthier-Villars et. fils, Paris, 1896, 254--256.


\bibitem {fujw1} \textsc{M. Fujiwara,} \textit{ $\ta$-congruent numbers}, in: Number Theory,
K. Gy\H{o}ry, A. Peth\H{o} and  V. S\'os (eds.), de Gruyter, 1997, pp. 235--241.

\bibitem {fujw2} \textsc{M. Fujiwara,} \textit{Some properties of  $\ta$-congruent numbers},
Natural Science Report, Ochanomizu University,  \textbf{118}(2) (2001) 1--8.


\bibitem{j-s} A. S. Janfada, S. Salami, {On $\ta$-congruent numbers on real quadratic number fields}, Kodai Math. J. \textbf{38} 352--364 (2015).


\bibitem {Halb-Hun} \textsc{Halbeisen, L. and Hungerb\"{u}hler, N.}, \textit{A Theorem of Fermat and Congruent Numbers Curves},  Hardy-Ramanujan Journal {\bf{41}} (2018) 15-21.


\bibitem {Hun} \textsc{Hungerb\"{u}hler, N.}, \textit{A proof of a conjecture of Lewis Carroll},
 Mathematics Magazine {\bf{69}} (1996), 182--184.

\bibitem {kan1} \textsc{M. Kan}, \textit{$\ta$-congruent numbers and elliptic curves}, Acta Arith. \textbf{94} (2000) 153--160.

\bibitem {kob1}\textsc {N. Koblitz}, \textit{Introduction to Elliptic curves and Modular Forms,}
Springer-Verlag, Graduate text in Mathematics 97, 2nd ed, Berlin, 1993.

\bibitem {mazur}\textsc {B. Mazur}, \textit{Modular curves and the Eisenstein ideal,} 
Pub. Math. IHES \textbf{47} (1978), 33--186.




\bibitem {slm1}{J. H. Silverman}, \emph{The Arithmetic of Elliptic Curves, second edition}
Graduate text in Mathematics, Springer-Verlag, Berlin, Vol. 106,  (2009).




\end{thebibliography}
\end{document}